\documentclass[12pt]{amsart}
\usepackage{amssymb}
\usepackage[all]{xy}
\usepackage{graphics}

\newtheorem{thm}{Theorem}[section]

\newtheorem{prop}[thm]{Theorem}

\newtheorem*{prob*}{Open problem}
\theoremstyle{definition}

\theoremstyle{remark}

\newtheorem*{rem*}{Remark}

\newcommand{\kringel}{\mathbin{\raise1pt\hbox{$\scriptstyle\circ$}}} 
\newcommand{\pkt}{\mathbin{\raise0pt\hbox{$\scriptstyle\bullet$}}}

\newcommand{\F}{\mathbb{F}}

\newcommand{\Q}{\mathbb{Q}}

\newcommand{\Z}{\mathbb{Z}}

\newcommand{\ov}{\overline}

\newcommand{\ra}{\rightarrow}  

\renewcommand{\phi}{\varphi}

\begin{document}


\title[]{A remark on an inequality for the prime counting function}
\author[D. Burde]{Dietrich Burde}
\date{\today}
\email{dietrich.burde@univie.ac.at}

\address{Fakult\"at f\"ur Mathematik\\
Universit\"at Wien\\
  Nordbergstrasse 15\\
  1090 Wien}

\begin{abstract}
We note that the inequalities $0.92 \frac{x}{\log(x)} <\pi(x)< 1.11 \frac{x}{\log(x)}$ do not
hold for all $x\ge 30$, contrary to some references. These estimates on $\pi(x)$ came up recently 
in papers on algebraic number theory.
\end{abstract}

\maketitle

\section{Chebyshev's estimates for $\pi(x)$}

Let $\pi(x)$ denote the number of primes not greater than $x$, i.e.,\\
\begin{align*}
\pi(x) & = \sum_{p\le x} 1. \\
\end{align*}
One of the first works on the function $\pi(x)$ is due to Chebyshev.
He proved (see \cite{CHE}) in $1852$ the following explicit inequalities for $\pi(x)$,
holding for all $x\ge x_0$ with some $x_0$ sufficiently large: \\
\begin{align*}
c_1 \frac{x}{\log(x)} & < \pi(x) < c_2 \frac{x}{\log(x)}, \\[0.3cm]
c_1 & = \log(2^{1/2}3^{1/3}5^{1/5}30^{-1/30})\approx 0.921292022934, \\[0.3cm]
 c_2 & =\frac{6}{5}c_1\approx 1.10555042752. \\
\end{align*}

This can be found in many books on analytic number theory (see for example \cite{APO}, \cite{CHA}, 
\cite{PRA} and \cite{TEN}).
But it seems that this result is sometimes cited incorrectly: it is claimed
that the estimates are valid for all $x\ge 30$. For example, in \cite{ELL}, page $21$ we read that\\
\begin{align*}
c_1 \frac{x}{\log(x)} & < \pi(x) < c_2 \frac{x}{\log(x)}, \quad \forall \, x\ge 30.\\
\end{align*}

But a quick numerical computation shows that this is wrong. To give an example, 
take $x=100$. Then we have $\pi(x)=25$ and 
\[
c_2 \frac{x}{\log(x)} \approx 24.00672250690558538515780234 <25.
\]
Actually, the inequality is far from true for small $x$. We have the following result:

\begin{prop}
Let $c_2\approx 1.10555042752$ be Chebyshev's constant.
Then the inequality \\
\begin{align*}
\pi(x) & < c_2 \frac{x}{\log(x)}
\end{align*}
is true for all $x\ge 96098$. For $x=96097$ it is false.
\end{prop}

\begin{proof}
In \cite{PAN} it is shown that\\
\begin{align*}
\pi(x) & <  \frac{x}{\log(x)-1.11},  \quad x\ge 4.\\
\end{align*}
The RHS is less or equal to $c_2x/\log(x)$ if and only if \\
\begin{align*}
x & \ge  \exp\left(\frac{1.11\cdot c_2}{c_2-1}\right)\approx 112005.18. \\
\end{align*}
This shows the claim for $x\ge 112006$. 
Since $x/\log(x)$ is a monotonously increasing
function it is enough to check the claimed estimate for intergers $x$ in the intervall
$[96098, 112006]$ by computer. For $x=96097$ we have $\pi(96097)=9260$ and
$c_2 x/\log (x)\approx 9259.92$.
\end {proof}

The incorrect inequality was also used in a former version of Khare's proof of Serre's
modularity conjecture for the level one case, see \cite{KHA}, \cite{KHA1}. Let $\F$ be a finite field
of characteristic $p$. The conjecture stated that
an odd, irreducible Galois representation $\rho\colon Gal(\ov{\Q}/\Q)\ra GL_2(\F)$ which is
unramified outside $p$ is associated to a modular form on $SL_2(\Z)$. Khare's proof is an
elaborate induction on $p$. Starting with a $p$ for which the conjecture is known one wants to
prove the conjecture for a larger prime $P$. Kahre's arguments do only work if $P$ and $p$ are
not Fermat primes, and if
\[
\frac{P}{p}\le a
\]
for certain values $a>1$, close to $1$. At this point Khare used the incorrect estimate on $\pi(x)$,
as explained above. Fortunately the proof easily could be repaired by using better estimates
on $\pi(x)$ provided by Rosser and Schoenfeld \cite{ROS}, and Dusart \cite{DU1}.\\
Indeed, P. Dusart proved inequalities for $\pi(x)$ which are much
better than Chebyshev's estimates. He verifies this for smaller $x$ numerically.
Nevertheless he claims in his thesis \cite{DU2}, that Chebyshev 
gave the following inequality\\
\begin{align*}
0.92 \frac{x}{\log(x)} & < \pi(x) < 1.11 \frac{x}{\log(x)}, \quad x\ge 30,\\
\end{align*}
which is equally wrong. \\
The question is: where lies the origin for this error ?
Chebyshev himself proved inequalities in \cite{CHE} with his constants $c_1$ and $c_2=\frac{6}{5}c_1$
indeed for all $x\ge 30$, but for inequalities involving $\psi(x)=\sum_{n\le x}\Lambda(n)$ instead of $\pi(x)$. 
His estimates concerning $\psi(x)$ seem to be correct for all $x\ge 30$. For example, he shows by
elementary means that, for all $x\ge 30$, \\
\begin{align*}
\psi(x) & < \frac{6}{5}c_1x+ \frac{5}{4\log(6)}\log^2(x)+\frac{5}{4}\log (x) +1, \\
\psi(x) & > c_1x-\frac{5}{2}\log (x)-1.\\
\end{align*}
To derive from this inequalities on $\pi(x)$ for $x\ge 30$, we have to estimate\\
\begin{align*}
\psi(x) & = \sum_{p\le x}\left[ \frac{\log(x)}{\log(p)} \right] \log(p).\\
\end{align*}
Using the estimates $[y]\le y<[y]+1\le 2[y]$ for $y\ge 1$ we obtain \\
\begin{align*}
\psi(x) & \le \pi(x)\log(x) \le 2\psi(x), \quad x\ge 2.\\
\end{align*}
On the RHS we cannot do easily much better than $2\psi(x)$. Hence we obtain\\
\begin{align*}
c_1 \frac{x}{\log(x)} & < \pi(x) < 2c_2 \frac{x}{\log(x)} , \quad x\ge 30.\\
\end{align*}
On the other hand we know that\\
\begin{align*}
\pi(x) & = \frac{\psi(x)}{\log(x)}+O\left( \frac{x}{\log^2(x)}\right), \quad x\ge 2,\\
\end{align*}
so that we obtain, as $x$ tends to infinity,\\
\begin{align*}
(c_1 +o(1))\frac{x}{\log(x)} & \le \pi(x) \le (c_2+o(1)) \frac{x}{\log(x)}.\\
\end{align*}
Chebyshev used these estimates to prove Bertrand's postulate: 
each interval $(n,2n]$ for $n\ge 1$ contains at least one prime.
Moreover his results were a first step towards the proof of the prime number theorem.

\section{Other estimates for $\pi(x)$}

There are many interesting inequalities on the function $\pi(x)$.
Let us first consider inequalities of the form 
\begin{align*}
A \frac{x}{\log(x)} & < \pi(x) < B \frac{x}{\log(x)}
\end{align*}
for all $x\ge x_0$, where $x_0$ depends on the constant $A\le 1$ and respectively on $B>1$.
On the LHS we can choose $A$ equal to $1$, if $x\ge 17$. In fact, we have \cite{DU2}
\begin{align*}
\frac{x}{\log(x)} & < \pi(x) ,  \quad \forall \; x\ge 17.
\end{align*}
Note that for $x=16.999$ we have $x/\log(x)\approx 6.0000257$, but $\pi(x)=6$.
Consider the RHS of the above inequalities: if we want to hold such inequalities on $\pi(x)$ 
for all $x\ge x_0$ with a smaller $x_0$, we need to enlarge the constant $B$. Conversely, if we need this inequality for
smaller $B$, we have to enlarge $x_0$. The prime number theorem ensures that we can choose $B$ as close to $1$
as we want, provided $x_0$ is sufficiently large.
The following result of Dusart \cite{DU1} enables us to derive adjusted versions for the above inequalities:
\begin{thm}[Dusart]\label{2.1}
For real $x$ we have the following sharp bounds:\\
\begin{align*}
\pi(x) & \ge \frac{x}{\log(x)}\left( 1+\frac{1}{\log(x)}+\frac{1.8}{\log^2(x)}\right), \quad x\ge 32299,\\[0.2cm]
\pi(x) & \le \frac{x}{\log(x)} \left( 1+ \frac{1}{\log(x)} +\frac{2.51}{\log ^2(x)}\right),
\quad x\ge 355991.\\
\end{align*}
\end{thm}
One can derive, for example, the following inequalities.\\
\begin{align*}
\pi(x) & < 1.095 \cdot \frac{x}{\log(x)},  \quad x\ge 284860,\\
\pi(x) & < 1.25506 \cdot \frac{x}{\log(x)},  \quad x\ge 17.\\
\end{align*}
Among other inequalities on $\pi(x)$ we mention the following ones:\\
\begin{align*}
\frac{x}{\log(x)-m} & < \pi(x) < \frac{x}{\log(x)-M} \\
\end{align*}
for all $x\ge x_0$ with real constants $m$ and $M$. They have been studied by various authors.
A good reference is the article \cite{PAN}. There it is shown, for example, that 
\begin{align*}
\pi(x) & > \frac{x}{\log(x)-\frac{28}{29}}, \quad x\ge 3299,\\[0.2cm]
\pi(x) & < \frac{x}{\log(x)-1.11}, \quad x\ge 4. 
\end{align*}
The second inequality can also be used to obtain results on our estimate $\pi(x)< B\frac{x}{\log(x)}$,
in particular for smaller $x$, where the second inequality of Theorem $\ref{2.1}$ is not valid.
However we have
\begin{align*}
 \frac{x}{\log(x)} \left( 1+ \frac{1}{\log(x)} +\frac{2.51}{\log ^2(x)}\right) & < \frac{x}{\log(x)-1.11}, 
\quad x\ge 28516. 
\end{align*}
For $x>10^6$ and $a=1.08366$ we can use \cite{PAN}
\[
\pi(x)<\frac{x}{\log (x)-a}.
\]
Here the upper bound of Dusart is better only as long as $x \ge 2846396$. \\
Finally we mention the book \cite{SMC}, providing many references on inequalities on $\pi(x)$,
and the recent article \cite{HAS}, where lower and upper bounds 
for $\pi(x)$ of the form $\frac{n}{H_n-c}$ are discussed, where $H_n=1+\frac{1}{2}+\cdots + \frac{1}{n}$.

\section{acknowledgement} 
We are grateful to the referee for drawing our attention to several
approximations of $\pi(x)$. We thank J. S\'andor for helpful remarks.

\end{document}